\documentstyle[12pt]{article}

\newlength{\abstractwidth}
\renewcommand{\thefootnote}{\fnsymbol{footnote}}
\renewcommand{\thanks}[1]{\footnote{#1}} 
\newcommand{\starttext}{ \setcounter{footnote}{0}
\renewcommand{\thefootnote}{\arabic{footnote}}}

\newcommand{\be}{\begin{equation}}
\newcommand{\bea}{\begin{eqnarray}}
\newcommand{\eea}{\end{eqnarray}} \newcommand{\ee}{\end{equation}}

\def\ba{\begin{eqnarray}}
\def\ea{\end{eqnarray}}

\def\o{\omega}

\def\det{{\rm det}}

\def\log{\,{\rm log}\,}

\def\o{\omega}

\def\e{\varepsilon}

\def\o{\omega}

\def\na{\nabla}

\def\R{{\bf R}}

\def\p{\prod}

\def\ddb{{\partial\bar\partial}}

\def\na{{\nabla}}

\def\[{{\bf [}}
\def\]{{\bf ]}}

\def\p{\partial}

\def\u{{\bf u}}
\def\a{{\bf a}}

\begin{document}
\starttext \baselineskip=15pt \setcounter{footnote}{0}
\newtheorem{theorem}{Theorem}
\newtheorem{lemma}{Lemma}
\newtheorem{definition}{Definition}
\newtheorem{proposition}{Proposition}
\newtheorem{corollary}{Corollary}

\begin{center}
{\Large \bf ESTIMATES FOR A GEOMETRIC FLOW FOR THE TYPE IIB STRING
\footnote{Work supported in part by the National Science Foundation Grants DMS-1855947 and DMS-1809582.}}
\bigskip\bigskip

\centerline{Teng Fei, Duong H. Phong, Sebastien Picard, and Xiangwen Zhang}

\medskip

\begin{abstract}

\smallskip

{\footnotesize
It is shown that bounds of all orders of derivative would follow from uniform bounds for the metric and the torsion $1$-form, for a flow in non-K\"ahler geometry which can be interpreted as either a flow for the Type IIB string or the Anomaly flow with source term and zero slope parameter. A key ingredient in the proof is a formulation of this flow unifying it with the Ricci flow, which was recently found.}
\end{abstract}

\end{center}

\section{Introduction}
\setcounter{equation}{0}

Anomaly flows were introduced in \cite{PPZ0} as a way of enforcing the conformally balanced condition in the Hull-Strominger system for supersymmetric compactifications of the heterotic string \cite{CHSW, H, S}. Their effectiveness can be inferred from their having produced an alternative and unified approach \cite{PPZ1} to the solutions found by Fu and Yau \cite{FY1, FY2}, which had originally required solving a complicated Monge-Amp\`ere type equation with gradients and very different estimates. Since then, they have also revealed themselves to be at the interface of many well-known questions in complex geometry: the special case with zero slope and no source can give another proof of the Calabi conjecture as well as a test of the K\"ahler property \cite{PPZ2}; it also turns out to be a natural generalization of the Ricci flow to the complex and non-K\"ahler setting \cite{FP}, which had been considered independently in \cite{ST, U}; and we shall see shortly below that the version with zero slope but with source terms provides a natural geometric flow for supersymmetric compactifications of the Type IIB string, which had been formulated in \cite{GMPT,T} and \cite{TY1, TY2, TY3}. Anomaly flows also fit naturally in the broad theme of curvature flows in non-K\"ahler geometry, to which also belong many flows of great current interest \cite{BV1, BV2, BX, L, LW,  ST, U}.

\medskip
While Anomaly flows have been worked out in a number of examples (see e.g. \cite{FHP, FP, FeY, PPZ1, PPZ3} and the recent case of nilmanifolds \cite{PU}), their understanding is still very incomplete. In particular, no general criterion for their long-time existence or convergence is as yet known. For this, it would be valuable to know what minimum set of estimates would imply estimates to all orders. This issue is of particular interest for flows in non-K\"ahler geometry, as the precise role of torsion terms is not yet fully understood. Even for scalar equations, $C^1$ estimates often stand apart and require different tools (see e.g. \cite{DK, PPZ4}).

\medskip
The goal of the present paper is to address this issue for the following flow. Let $X$ be a compact $n$-dimensional complex manifold, equipped with a nowhere vanishing holomorphic $(n,0)$-form $\Omega$, and a conformally balanced Hermitian metric $\o_0$, i.e., $d(\|\Omega\|_{\o_0}\omega_0^{n-1})=0$ where $\|\Omega\|_{\o_0}$ is the norm of $\Omega$ with respect to $\o_0$, that is $i^{n^2}\Omega\wedge\bar\Omega=\|\Omega\|_\o^2(\o^n/n!)$. Let $\Psi$ be a smooth given real and closed $(n-1,n-1)$-form. We consider the flow $t\to\o(t)$ defined by
\bea
\label{TypeIIB}
\p_t(\|\Omega\|_\o\o^{n-1})=i\p\bar\p\o^{n-2}-\Psi.
\eea
A fundamental property of this flow is that it preserves the conformally balanced condition $d(\|\Omega\|_\o\o^{n-1})=0$.
If we specialize to $n=3$, we recognize this flow as the flow of the Hermitian metric in the Anomaly flow defined in \cite{PPZ0}, with slope $\alpha'
=0$ and source $\Psi$. Also, recall that the equations for supersymmetric compactifications for the Type IIB string have been worked out in \cite{GMPT} and \cite{T}. In the simplest form studied in \cite{TY1,TY2,TY3}, they can be expressed as follows. In this case, one sets again $n=3$, and looks for a Hermitian metric $\hat\o$ on $X$ satisfying the following equations
\bea
d\hat\omega^2=0,
\qquad
&&
i\p\bar\p  (\|\Omega\|_{\hat \o}^{-2}\hat\o)={1\over 2}\rho_B
\eea
where $\rho_B$ is the Poincar\'e dual to a given linear combination of holomorphic $2$-cycles. If we set
$\o=\|\Omega\|_{\hat\o}^{-2}\hat\o$,
then $\|\Omega\|_{\o}=\|\Omega\|_{\hat\o}^4$, then the preceding equations can be rewritten as
\bea
d(\|\Omega\|_{\o}\omega^2)=0,
\qquad
i\p\bar\p\o={1\over 2}\rho_B.
\eea
This just means that $\o$ is a stationary point of the flow (\ref{TypeIIB}) with source $\Psi={1\over 2}\rho_B$, and the flow provides a natural parabolic approach to solving the equations of the Type IIB compactifications (see also \cite{GF}). Henceforth, we shall refer to (\ref{TypeIIB}) flow as the Type IIB flow, reserving the name ``Anomaly flow" for the cases involving a non-zero parameter $\alpha'$.

\medskip
To state precisely our main results, we describe our conventions. If we write $\omega = i g_{\bar{k} j} dz^j \wedge d \bar{z}^k$, the torsion of $\omega$ is
\be
T^m{}_{jp} = g^{m \bar{q}} (\partial_j g_{\bar{q} p} - \partial_p g_{\bar{q} j}).
\ee
The torsion 1-form $\tau$ has components
\be
\tau_i = T^p{}_{pi}.
\ee
It is well-known (e.g. \cite{PPZ5}) that for conformally balanced metrics, the identity
\be
\tau_i = \partial_i \log \| \Omega \|_\omega
\ee
holds. Our main theorem is then

\begin{theorem} \label{thm-long-time}
Let $(X,\omega_0)$ be a compact Hermitian manifold with nowhere vanishing holomorphic $(n,0)$ form $\Omega$ satisfying $d ( \| \Omega \|_{\omega_0} \omega_0^{n-1})=0$ and given closed $\Psi \in \Lambda^{n-1,n-1}(X)$. Let $\omega(t)$ evolve by the Type IIB flow
\be\label{af-s1}
\p_t \left( \|\Omega\|_{\o} \o^{n-1} \right) = i\ddb \omega^{n-2} - \Psi, \ \ \ \omega(0)=\omega_0.
\ee
Suppose on $[0,T_0]$, we have the estimate
\be
K_1^{-1} \omega_0 \leq \omega(t) \leq K_1 \omega_0, \ \ \ |\tau| \leq K_2
\ee
for $K_1,K_2>0$. Then, there exist constants $C>0$ and $0< \alpha<1$ depending on $K_1, K_2$, $(X, \omega_0,\Omega)$ and $\Psi$ such that
\be\label{high-order-est}
\| g \|_{C^{2+\alpha, 1+ \alpha/2}(X \times (0, T_0])} \leq C.
\ee
\end{theorem}

As an intermediate consequence of the above theorem, we can provide a criterion for the long-time existence of the flow (\ref{af-s1})

\begin{corollary}
Under the same hypotheses as in the above theorem, the flow can be extended to $[0,T_0+\e)$ for some $\e>0$.
\end{corollary}

In the case with no source, we can weaken our assumptions:

\begin{theorem} \label{thm-zero-source}
Let $(X,\omega_0)$ be a compact Hermitian manifold with nowhere vanishing holomorphic $(n,0)$ form $\Omega$ satisfying $d ( \| \Omega \|_{\omega_0} \omega_0^{n-1})=0$. Let $\omega(t)$ evolve by the Anomaly flow
\be
\p_t \left( \|\Omega\|_{\o} \o^{n-1} \right) = i\ddb \omega^{n-2}, \ \ \ \omega(0)=\omega_0.
\ee
Suppose on $[0,T_0]$, we have the estimate
\be
\omega(t) \leq K_1 \omega_0, \ \ \ |\tau| \leq K_2
\ee
for $K_1,K_2>0$. Then, there exist constants $C>0$ and $0< \alpha<1$ depending on $K_1, K_2$, $(X, \omega_0,\Omega)$ and $\Psi$ such that
\be\label{high-order-est}
\| g \|_{C^{2+\alpha, 1+ \alpha/2}(X \times (0, T_0])} \leq C.
\ee
Moreover, there exists $\epsilon>0$ such that the flow can be extended to $[0,T_0+\epsilon)$.
\end{theorem}

We stress several aspects of the estimates obtained in the above theorems. First, it is remarkable that only a bound for $\tau$ was needed, and not for the whole torsion tensor $T$. Second, unlike the Shi-type estimates obtained in \cite{PPZ5} assuming bounds for $|Rm|_\o+|T|_\o+|DT|_\o$, they are independent of the time interval $T_0$. Third, they allow for a non-trivial source term $\Psi$. The incorporation of sources had turned out to be surprisingly difficult in the approach of \cite{PPZ5}, and a key catalyst for the approach in the present paper is the new formulation of the Type IIB flow found in \cite{FP} which unifies it with the Ricci flow. This allows in particular the application of Calabi-type identities for the connection, as introduced in \cite{PSS}. Once the $C^1$ estimates for the metric have been obtained, we can apply the general Schauder theory for systems of quasilinear parabolic system to obtain all the higher order estimates \cite{LU}.

\section{First Order Estimate of the Metric}
\setcounter{equation}{0}

The bulk of the paper is devoted to the proof of Theorem \ref{thm-long-time}.

\bigskip

\par To prove Theorem \ref{thm-long-time}, following \cite{FP}, we consider the flow of
\be
\eta_{\bar{k} j} = \| \Omega \|_\omega g_{\bar{k} j}.
\ee
The form $\eta = i \eta_{\bar{k} j} dz^j \wedge d \bar{z}^k$ satisfies
\be
d( \| \Omega \|^2_\eta \eta^{n-1}) = 0.
\ee
In \cite{FP}, it is shown that $\eta$ evolves by
\be\label{uf}
\p_t \eta_{\bar k j} = - \tilde{R}_{\bar k j}(\eta) - {1\over 2} T_{\bar k p q}(\eta) \, \bar{T}_{j}{}^{pq}(\eta) - \Phi_{\bar k j}
\ee
where $\tilde{R}_{\bar k j}(\eta)$ is the second Chern-Ricci curvature of the metric $\eta$ and $\Phi(z,\eta_{\bar{k} j}) \in \Lambda^{1,1}(X)$ involves combinations of the given $\Psi \in \Lambda^{n-1,n-1}(X)$ and $\eta_{\bar{k} j}$.
\medskip
\par This evolution equation for $\eta$ is simpler than the one for $\omega$, and we will use it to obtain estimates on $\eta$ and then deduce estimates on $\omega$. For this, we note a relation between their torsion 1-forms. Taking the determinant of the defining relation, we see that
\be
\| \Omega \|^2_\eta = \| \Omega \|_g^{2-n}.
\ee
Therefore
\be
\partial_i \log \| \Omega \|_\eta = \partial_i \log \| \Omega \|_g^{1-(n/2)} = {2-n \over 2} \tau_i.
\ee
It follows that if
\be
K^{-1} \omega_0 \leq \omega(t) \leq K \omega_0, \ \ \ |\tau| \leq K
\ee
then
\be \label{eta-assumption}
C^{-1} \eta_0 \leq \eta(t) \leq C \eta_0, \ \ \ |\nabla \log \| \Omega \|_\eta|_\eta \leq C,
\ee
and for the remainder of the paper, we will assume (\ref{eta-assumption}) and work exclusively with the metrics $\eta(t)$. Note that if $\eta$ is bounded in the $C^k$ norm and $\eta \geq C^{-1} \eta_0$, then
\be
g_{\bar{k} j} = \| \Omega \|_\eta^{2/(n-2)} \eta_{\bar{k} j}
\ee
is also bounded in the $C^k$ norm.

\medskip

In the remaining part of this section, we follow the calculation in \cite{PSS} (see also \cite{PSSW, ZZ}) to derive a gradient estimate for the metric. Let $(X, g)$ be a compact Hermitian manifold. We will work with the flow
\be \label{typeiib-flow}
\p_t g_{\bar k j} = - \tilde{R}_{\bar k j} - {1\over 2} T_{\bar k p q} \, \bar{T}_{j}{}^{pq} - \Phi(z,g(t))_{\bar k j}.
\ee
In comparison to the flow (\ref{uf}), here we write $g_{\bar{k} j}$ instead of $\eta_{\bar{k} j}$. For abbreviation, we will denote $\hat g= g(0) = g_0$ and $g= g(t)$. Let
\bea
h^\alpha{}_{\beta}= \hat g^{\alpha \bar\gamma} \, g_{\bar \gamma \beta}
\eea
In the following computation, we will use $\nabla$ and $\hat\nabla$ to denote the Chern connections, $\theta$ and $\hat\theta$ to denote the connection 1-forms, $R$ and $\hat{R}$ to denote the curvatures, with respect to the metrics $g$ and $\hat{g}$. We write
\be
\Delta = g^{p \bar{q}} \nabla_p \nabla_{\bar{q}},
\ee
for the Laplacian, and
\be \label{Lambda-defn}
i \Lambda_\omega \Phi = g^{j \bar{k}} \Phi_{\bar{k} j}, \quad \omega = i g_{\bar{k} j} \, dz^j \wedge d \bar{z}^k
\ee
for the contraction operator. Our curvature conventions are
\be \label{chern-curv}
R_{\bar{k} j}{}^p{}_q = - \partial_{\bar{k}} \Gamma^p_{jq}, \quad \Gamma^p_{jq} = g^{p \bar{\alpha}} \partial_j g_{\bar{\alpha} q},
\ee
and
\be
\tilde{R}_{\bar{p} q} = g_{\bar{p} \ell} g^{j \bar{k}} R_{\bar{k} j}{}^\ell{}_q
\ee
for the second Chern-Ricci curvature.

\subsection{Evolution of relative endomorphism}
In this section, we derive the following evolution equation.

\begin{proposition}
Along the flow (\ref{typeiib-flow}), we have
\be
(\p_t - \Delta) {\rm Tr} \, h = -  g^{q \bar{p}} h^{-1}{}^\gamma{}_\mu \hat{\nabla}_q h^\mu{}_j \hat{\nabla}_{\bar{p}} h^j{}_\gamma   - g^{q \bar{p}} \hat{R}_{\bar{p} q}{}^\alpha{}_j h^j{}_\alpha    - {1 \over 2} \hat{g}^{j \bar{k}} T_{\bar{k} pq} \bar{T}_j{}^{pq} - i \Lambda_{\hat{\omega}} \Phi.
\ee
\end{proposition}

\par \noindent {\it Proof:} By definition,
\be
\p_t {\rm Tr} \, h = \hat{g}^{p \bar{q}} \p_t g_{\bar{q} p},
\ee
which by the equation for the flow is
\be \label{trh-evol1}
\p_t {\rm Tr} \, h = - \hat{g}^{j \bar{k}} \tilde{R}_{\bar{k} j} - {1 \over 2} \hat{g}^{j \bar{k}} T_{\bar{k} pq} \bar{T}_j{}^{pq} - i \Lambda_{\hat{\omega}} \Phi.
\ee
In general, we can compute
\bea \label{difference-connections}
\theta
&=& g^{-1}\, \p g
= h^{-1} \, \hat{g}^{-1} \, \p(\hat{g} \, h)
= h^{-1} \, \p h + h^{-1}\, \hat\theta \, h\\
\nonumber
&=&h^{-1}\, \p h  + h^{-1}\, \hat\theta \, h -  h^{-1}\, h \,\hat\theta + \hat\theta\\
\nonumber
&=& \hat\theta +h^{-1}\, \hat\nabla h.
\eea
The definition of the curvature of the Chern connection (\ref{chern-curv}) implies
\be
R_{\bar{p} q}{}^\alpha{}_\beta = \hat{R}_{\bar{p} q}{}^\alpha{}_\beta - \partial_{\bar{p}} ( h^{-1} \hat{\nabla}_q h)^\alpha{}_\beta
\ee
and
\be
\tilde{R}_{\bar{k} j} = g_{\bar{k} \alpha} g^{q \bar{p}} \hat{R}_{\bar{p} q}{}^\alpha{}_j - g_{\bar{k} \alpha} g^{q \bar{p}} \hat{\nabla}_{\bar{p}} ( h^{-1} \hat{\nabla}_q h)^\alpha{}_j.
\ee
Substituting this identity into (\ref{trh-evol1}), we obtain
\be
\p_t {\rm Tr} \, h =  g^{q \bar{p}} \hat{\nabla}_{\bar{p}} ( h^{-1} \hat{\nabla}_q h)^\alpha{}_j h^j{}_\alpha - g^{q \bar{p}} \hat{R}_{\bar{p} q}{}^\alpha{}_j h^j{}_\alpha    - {1 \over 2} \hat{g}^{j \bar{k}} T_{\bar{k} pq} \bar{T}_j{}^{pq} - i \Lambda_{\hat{\omega}} \Phi.
\ee
Expanding the first term gives
\bea
\p_t {\rm Tr} \, h &=& g^{q \bar{p}} (h^{-1})^\alpha{}_\mu \hat{\nabla}_{\bar{p}} \hat{\nabla}_q h^\mu{}_j h^j{}_\alpha  +  g^{q \bar{p}} (\hat{\nabla}_{\bar{p}} h^{-1})^\alpha{}_\mu (\hat{\nabla}_q h)^\mu{}_j h^j{}_\alpha  \nonumber\\
&&- g^{q \bar{p}} \hat{R}_{\bar{p} q}{}^\alpha{}_j h^j{}_\alpha    - {1 \over 2} \hat{g}^{j \bar{k}} T_{\bar{k} pq} \bar{T}_j{}^{pq} - i \Lambda_{\hat{\omega}} \Phi.
\eea
Since
\be
\hat{\nabla}_{\bar{p}} (h^{-1}){}^\alpha{}_\mu = - (h^{-1}){}^\alpha{}_\nu \hat{\nabla}_{\bar{p}} h^\nu{}_\gamma (h^{-1}){}^\gamma{}_\mu,
\ee
we obtain
\bea
\p_t {\rm Tr} \, h &=& g^{q \bar{p}} \partial_p \partial_{\bar{q}} {\rm Tr} \,  h -  g^{q \bar{p}} (h^{-1}){}^\gamma{}_\mu \hat{\nabla}_q h^\mu{}_j \hat{\nabla}_{\bar{p}} h^j{}_\gamma   \nonumber\\
&&- g^{q \bar{p}} \hat{R}_{\bar{p} q}{}^\alpha{}_j h^j{}_\alpha    - {1 \over 2} \hat{g}^{j \bar{k}} T_{\bar{k} pq} \bar{T}_j{}^{pq} - i \Lambda_{\hat{\omega}} \Phi.
\eea
Here we used that
\be
\hat{\nabla}_{\bar{p}} \hat{\nabla}_q h^i{}_i = \partial_{\bar{p}} (\partial_q h^i{}_i + \hat{\Gamma}^i{}_{qj} h^j{}_i - h^i{}_j \hat{\Gamma}^j{}_{qi} )= \partial_{\bar{p}} \partial_q {\rm Tr} \, h.
\ee
This proves the identity. Q.E.D.
\bigskip
\par This identity allows us to use the quantity ${\rm Tr} \, h$ to produce a negative term involving the $C^1$ norm of the metric.

\begin{proposition}
Let $g(t)$ evolve by Type IIB flow (\ref{typeiib-flow}). Suppose on $[0,T]$, we have the estimate
\be
K^{-1} \hat{g} \leq g(t) \leq K \hat{g},
\ee
for $K>0$. Then
\be  \label{Trh-est}
(\p_t - \Delta) {\rm Tr} \, h \leq - {1 \over K^2} |\hat{\nabla} g|^2 + C,
\ee
on $[0,T]$, where $C$ depends on $K$, $(X,\hat{g})$, $\Phi$.
\end{proposition}

{\it Proof:} We can estimate
\be
g^{q \bar{p}} (h^{-1}){}^\gamma{}_\mu \hat{\nabla}_q h^\mu{}_j \hat{\nabla}_{\bar{p}} h^j{}_\gamma = g^{p \bar{q}} g^{i \bar{k}} \hat{g}^{\ell \bar{k}} \hat{\nabla}_q g_{\bar{k} \ell}  \hat{\nabla}_q g_{\bar{k} i} \geq K^{-2} |\hat{\nabla} g|^2.
\ee
We also note that
\be
- {1 \over 2} \hat{g}^{j \bar{k}} T_{\bar{k} pq} \bar{T}_j{}^{pq} \leq 0.
\ee
Altogether, we obtain
\be
(\p_t - \Delta) {\rm Tr} \, h \leq -K^{-2} |\hat{\nabla} g|^2 - g^{q \bar{p}} \hat{R}_{\bar{p} q}{}^\alpha{}_j h^j{}_\alpha - i \Lambda_{\hat{\omega}} \Phi,
\ee
which implies the estimate. Q.E.D.

\subsection{$C^1$ estimate for the metric}

As in \cite{Yau} and \cite{PSS}, we consider
\bea
S=|\nabla h \, h^{-1} |^2_{g}= g^{m\bar \gamma} \, g_{\bar \mu \beta} \, g^{\ell \bar \alpha} \, (\nabla_m h \, h^{-1} )^\beta{}_{\ell} \, \overline{(\nabla_\gamma h  \, h^{-1} )^\mu{}_{\alpha}}.
\eea
The guiding principle of the computations below is that $\na h h^{-1}$ is essentially a connection.
Note that by the identity (which can be derived by a computation similar to (\ref{difference-connections}))
\be
\hat{\theta} = \theta - \nabla h h^{-1},
\ee
we can write
\be
S = |\theta - \hat{\theta}|^2_g.
\ee
By the identity
\be
\Gamma^k_{ij} - \hat{\Gamma}^k_{ij} =  g^{k \bar{\ell}} \partial_i g_{\bar{\ell} j} - \hat{\Gamma}^k_{ij} = g^{k \bar{\ell}} \hat{\nabla}_i g_{\bar{\ell} j},
\ee
we can also write
\be
S = |\hat{\nabla} g|^2_g.
\ee
In this section, we will show the following estimate.

\begin{proposition} \label{prop-S-estimate}
Let $g(t)$ evolve by Type IIB flow (\ref{typeiib-flow}). Suppose on $[0,t_0]$, we have the estimate
\be
K^{-1} \hat{g} \leq g(t) \leq K \hat{g}, \ \ |T|^2 \leq K
\ee
for $K>0$. Then
\be
(\partial_t - \Delta)S \leq C S + C,
\ee
on $[0,t_0]$, where $C$ depends on $K$, $(X,\hat{g})$, $\Phi$.
\end{proposition}

Given Proposition \ref{prop-S-estimate} and (\ref{Trh-est}), we conclude
\begin{theorem} \label{thm-c1-estimate}
Let $g(t)$ evolve by Type IIB flow (\ref{typeiib-flow}). Suppose on $[0,t_0]$, we have the estimate
\be
K^{-1} \hat{g} \leq g(t) \leq K \hat{g}, \ \ |T|^2 \leq K
\ee
for $K>0$. Then
\be
| \hat{\nabla} g|^2_g \leq C
\ee
on $[0,t_0]$, where $C$ depends on $K$, $(X,\hat{g})$, $\Phi$, $g(0)$.
\end{theorem}

{\it Proof:} Let
\be
G = S- A {\rm Tr} \, h,
\ee
where $A \gg 1$ will be chosen depending on $K$, $(X,\hat{g})$, $\Phi$. Combining Proposition \ref{prop-S-estimate} and (\ref{Trh-est}), we obtain
\be
(\p_t - \Delta) G \leq C_1 S + C - {A \over K^2} S + AC.
\ee
Choosing $A = K^2 C_1 + K^2$,
\be
(\p_t - \Delta)G \leq AC - S.
\ee
Let $(p,t) \in X \times [0,T]$, be a point where $G$ attains its maximum with $t>0$. Then
\be
0 \leq AC - S(p,t),
\ee
which implies $S(p,t) \leq AC$, and hence
\be
S = G + A {\rm Tr} \, h \leq G(p,t) + C \leq C.
\ee
Since $S = |\hat{\nabla} g|^2$, this proves the desired estimate. Q.E.D.

\subsubsection{General formula for the evolution of $S$}
From \cite{PSS}, under any flow, we have the general formula
\bea \label{evol-S}
(\p_t - \Delta) S &=&
-|\bar\nabla(\nabla h \, h^{-1} )|^2-|\nabla(\nabla h \, h^{-1} )|^2\\\nonumber
&&
+ g^{m\bar \gamma} \left\langle\, (\p_t - \Delta) (\nabla_m h \, h^{-1}  ), \, \nabla_{\gamma} h \, h^{-1}  \, \right\rangle
+ g^{m\bar\gamma} \left\langle \, \nabla_m h \, h^{-1} , \, (\p_t - \Delta ) (\nabla_\gamma h \, h^{-1} ) \, \right\rangle \\
\nonumber
&&
- (\nabla_m h \, h^{-1} )^\beta{}_{\ell}\, \overline{(\nabla_\gamma h \, h^{-1} )^\mu{}_{\alpha}} \cdot
\bigg\{
\left(h^{-1} \, \dot{h} + \tilde{R}\right)^{m\bar \gamma} g_{\bar \mu \beta} \, g^{\ell\bar\alpha}\\
\nonumber
&&
 - g^{m\bar \gamma} \left(h^{-1} \, \dot{h}+ \tilde{R}\right)_{\bar \mu \beta}\, g^{\ell \bar \alpha}+ g^{m\bar \gamma} \, g_{\bar \mu \beta} \, \left(h^{-1}\, \dot{h} + \tilde{R}\right)^{\ell \bar \alpha}
\bigg\}
\eea
We provide the derivation for completeness. We have
\be
\nabla_{\bar{q}} S = \langle \nabla_{\bar{q}} \nabla h h^{-1}, \nabla h h^{-1} \rangle + \langle \nabla h h^{-1}, \nabla_q \nabla h h^{-1} \rangle,
\ee
and
\bea
\Delta S &=& g^{p \bar{q}}  \langle \nabla_p \nabla_{\bar{q}} \nabla h h^{-1}, \nabla h h^{-1} \rangle + g^{p \bar{q}} \langle \nabla_{\bar{q}} \nabla h h^{-1}, \nabla_{\bar{p}} \nabla h h^{-1} \rangle  \nonumber\\
&&+ g^{p \bar{q}} \langle \nabla_p \nabla h h^{-1}, \nabla_q \nabla h h^{-1} \rangle + g^{p \bar{q}} \langle \nabla h h^{-1}, \nabla_{\bar{p}} \nabla_q \nabla h h^{-1} \rangle .
\eea
Therefore
\bea
\Delta S &=& | \bar{\nabla} (\nabla h h^{-1})|^2 + |\nabla (\nabla h h^{-1})|^2 + \langle \Delta (\nabla h h^{-1}), \nabla h h^{-1} \rangle \nonumber\\
&&+ g^{p \bar{q}} \langle \nabla h h^{-1}, \nabla_{\bar{p}} \nabla_q (\nabla h h^{-1}) \rangle.
\eea
Our curvature convention is
\be
[\nabla_j,\nabla_{\bar{k}}]W_i = - R_{\bar{k} j}{}^p{}_i W_p, \ \ \ [\nabla_j,\nabla_{\bar{k}}]W_{\bar{i}} = R_{\bar{k} j \bar{i}}{}^{\bar{p}} W_{\bar{p}}.
\ee
Therefore
\bea
\nabla_{\bar{p}} \nabla_q (\nabla_i h h^{-1})^\alpha{}_\beta &=& \nabla_q \nabla_{\bar{p}}  (\nabla_i h h^{-1})^\alpha{}_\beta + R_{\bar{p} q}{}^j{}_i (\nabla_j h h^{-1})^\alpha{}_\beta - R_{\bar{p} q}{}^\alpha{}_\gamma (\nabla_i h h^{-1})^\gamma{}_\beta \nonumber\\
&&+ R_{\bar{p} q}{}^\gamma{}_\beta (\nabla_i h h^{-1})^\alpha{}_\gamma.
\eea
It follows that
\bea
\Delta S &=& | \bar{\nabla} (\nabla h h^{-1})|^2 + |\nabla (\nabla h h^{-1})|^2 + \langle \Delta (\nabla h h^{-1}), \nabla h h^{-1} \rangle \nonumber\\
&&+ \langle \nabla h h^{-1}, \Delta (\nabla h h^{-1}) \rangle +  g^{k \bar{i}} \langle \nabla_k h h^{-1} ,\tilde{R}^{j}{}_i (\nabla_j h h^{-1})^\alpha{}_\beta \rangle \nonumber\\
&&- g^{k \bar{i}} \langle \nabla_k h h^{-1}, \tilde{R}^\alpha{}_\gamma (\nabla_i h h^{-1})^\gamma{}_\beta \rangle\nonumber\\
&&+ g^{k \bar{i}} \langle \nabla_k h h^{-1}, \tilde{R}^\gamma{}_\beta (\nabla_i h h^{-1})^\alpha{}_\gamma \rangle.
\eea
On the other hand,
\be
\p_t S =  \p_t \bigg[ g^{k \bar{i}} g_{\bar{\mu} \alpha} g^{\beta \bar{\nu}} (\nabla_k h h^{-1})^\alpha{}_\beta \overline{(\nabla_i h h^{-1})^{\mu}{}_\nu} \bigg].
\ee
We note
\be
\dot{h}^\alpha{}_\beta = \hat{g}^{\alpha \bar{\mu}} \dot{g}_{\bar{\mu} \beta},  \ \ (h^{-1} \dot{h})^\alpha{}_\beta = g^{\alpha \bar{\mu}} \dot{g}_{\bar{\mu} \beta}.
\ee
Therefore
\bea
\p_t S &=& \langle \p_t (\nabla h h^{-1}), \nabla h h^{-1} \rangle + \langle \nabla h h^{-1}, \p_t (\nabla h h^{-1}) \rangle \nonumber\\
&& -  (h^{-1} \dot{h} )^k{}_\ell g^{\ell \bar{i}}  g_{\bar{\mu} \alpha} g^{\beta \bar{\nu}} (\nabla_k h h^{-1})^\alpha{}_\beta \overline{(\nabla_i h h^{-1})^{\mu}{}_\nu} \nonumber\\
&& + g^{k \bar{i}} g_{\bar{\mu} p} (h^{-1} \dot{h})^{p}{}_\alpha g^{\beta \bar{\nu}} (\nabla_k h h^{-1})^\alpha{}_\beta \overline{(\nabla_i h h^{-1})^{\mu}{}_\nu} \nonumber\\
&&- g^{k \bar{i}} g_{\bar{\mu} \alpha} (h^{-1} \dot{h})^\beta{}_p g^{p \bar{\nu}} (\nabla_k h h^{-1})^\alpha{}_\beta \overline{(\nabla_i h h^{-1})^{\mu}{}_\nu}
\eea
Putting this together, we obtain (\ref{evol-S}).

\subsubsection{Evolution of $S$ along Type IIB flow}
In this section, we derive the expression for the evolution of $S$ along the flow (\ref{typeiib-flow}). For this, we start by using
\be
(h^{-1} \, \dot{h})^\gamma{}_{\beta} = g^{\gamma \bar{\nu}} \dot{g}_{\bar{\nu} \beta},
\ee
to rewrite the flow as
\bea \label{evol-h}
(h^{-1} \, \dot{h})^\gamma{}_{\beta}+ \tilde{R}^\gamma{}_{\beta}= Q^\gamma{}_{\beta} -  \Phi^\gamma{}_{\beta}.
\eea
Here $Q_{\bar \gamma \beta} = - {1\over 2} T_{\bar \gamma pq} \, \bar{T}_\beta{}^{pq}$. Next, we compute $\p_t (\nabla_m h\, h^{-1})^{\alpha}{}_{\beta}$. Similarly as (\ref{difference-connections}), we have
\be \label{diff-connections2}
\hat\theta = \theta - \nabla h\, h^{-1}.
\ee
Therefore
\be
\p_t (\nabla_m h\, h^{-1})^{\alpha}{}_{\beta} = \p_t (\Gamma^\alpha_{m \beta} - \hat{\Gamma}^\alpha_{m \beta}) = \p_t (g^{\alpha \bar{\gamma}} \partial_m g_{\bar{\gamma} \beta}).
\ee
It follows that
\be
\p_t (\nabla_m h\, h^{-1})^{\alpha}{}_{\beta} = g^{\alpha \bar{\gamma}} (\partial_m \dot{g}_{\bar{\gamma} \beta} - \dot{g}_{\bar{\gamma}q} \Gamma^q_{m\beta}) = g^{\alpha \bar{\gamma}} \nabla_m \dot{g}_{\bar{\gamma} \beta}.
\ee
Thus
\bea \label{evol-nabla-g1}
\p_t (\nabla_m h\, h^{-1})^{\alpha}{}_{\beta}
&=& \nabla_m (h^{-1}\, \dot h)^{\alpha}{}_\beta
= g^{\alpha\bar\gamma} \nabla_m \left(- \tilde{R}_{\bar \gamma \beta} + Q_{\bar \gamma \beta} - \Phi_{\bar \gamma \beta} \right) \\\nonumber
&=& -\nabla_m \tilde{R}^\alpha{}_{\beta} + \nabla_m Q^\alpha{}_\beta - \nabla_m  \Phi^\alpha{}_\beta.
\eea

We need to compare this expression with $\Delta (\nabla_m h \, h^{-1})$. Using the relation between connection 1-forms (\ref{diff-connections2}), we obtain
\bea
\hat{R}_{\bar q m}{}^\alpha{}_\beta = R_{\bar q m}{}^\alpha{}_\beta +  \bar \partial_{q} \left(\nabla_m h \, h^{-1} \right)^\alpha{}_\beta.
\eea
It follows that
\bea
\Delta( \nabla_m h\, h^{-1})^{\alpha}{}_{\beta} =
g^{p\bar q} \nabla_p \nabla_{\bar q} ( \nabla_m h\, h^{-1})^{\alpha}{}_{\beta}  = \nabla^{\bar q} \, \hat{R}_{\bar q m}{}^\alpha{}_\beta - \nabla^{\bar q} \, R_{\bar q m}{}^\alpha{}_{\beta}.
\eea
Recall that, for the Chern connections of general Hermitian metrics, we have the following first Bianchi identies
\begin{eqnarray}
\label{1st}
R_{\bar \ell m \bar k j} = R_{\bar \ell j \bar k m} + \nabla_{\bar \ell} \, T_{\bar k jm}, \ \ R_{\bar \ell m \bar k j} = R_{\bar k m \bar \ell j} + \nabla_m \bar{T}_{j\bar k \bar \ell}
\end{eqnarray}
and also the second Bianchi identities
\begin{eqnarray}
\label{2nd-unbar}
&&
\nabla_m R_{\bar k j}{}^\alpha{}_{\beta} = \nabla_j R_{\bar k m}{}^\alpha{}_{\beta} + T^r{}_{jm} \, R_{\bar k r}{}^\alpha{}_{\beta}, \ \ \nabla_m R_{\bar k j \bar p q} = \nabla_j R_{\bar k m \bar p q} + T^r{}_{jm} \, R_{\bar k r \bar pq}\\
\label{2nd-bar}
&&
\nabla_{\bar m} R_{\bar k j }{}^{\alpha}{}_{\beta} = \nabla_{\bar k}R_{\bar m j}{}^\alpha{}{}_{\beta}+ \bar{T}^{\bar r}{}_{\bar k \bar m} \, R_{\bar r j}{}^\alpha{}_{\beta}, \ \ \nabla_{\bar m} R_{\bar k j \bar p q}= \nabla_{\bar k }R_{\bar m j \bar p q} + \bar{T}^{\bar r}{}_{\bar k \bar m} \, R_{\bar r j \bar p q}.
\end{eqnarray}
Using (\ref{2nd-unbar}), we can compute
\bea
\nabla^{\bar q} R_{\bar q m}{}^\alpha{}_{\beta}&=&
g^{p\bar q}\nabla_p R_{\bar q m}{}^\alpha{}_{\beta}
= g^{p\bar q}(\nabla_m R_{\bar q p}{}^\alpha{}_{\beta} + T^r{}_{mp}\, R_{\bar q r}{}^\alpha{}_\beta)\\\nonumber
&=&
\nabla_m \tilde{R}^\alpha{}_\beta + g^{p\bar q} T^r{}_{mp}\, R_{\bar q r}{}^\alpha{}_\beta
\eea
It follows that
\bea
\Delta( \nabla_m h\, h^{-1})^{\alpha}{}_{\beta} = -\nabla_m \tilde{R}^\alpha{}_\beta +
\nabla^{\bar q} \, \hat{R}_{\bar q m}{}^\alpha{}_\beta- g^{p\bar q} T^r{}_{mp}\, R_{\bar q r}{}^\alpha{}_\beta.
\eea
Thus, combining this with (\ref{evol-nabla-g1}), we obtain
\bea
(\p_t - \Delta)( \nabla_m h\, h^{-1})^{\alpha}{}_{\beta}
=
-\nabla^{\bar q} \, \hat{R}_{\bar q m}{}^\alpha{}_\beta
+g^{p\bar q} T^r{}_{mp}\, R_{\bar q r}{}^\alpha{}_\beta
+\nabla_m Q^\alpha{}_\beta - \nabla_m \Phi^\alpha{}_\beta
.
\eea
Substituting this expression and (\ref{evol-h}) into the general formula (\ref{evol-S}) for $S$, we obtain
\bea
(\p_t - \Delta) S &=&
-|\bar\nabla(\nabla h \, h^{-1} )|^2 -|\nabla(\nabla h \, h^{-1} )|^2\\\nonumber
&& + (I) +(II) + (III)+ (IV) + (V)+(VI)
\eea
where
\bea
(I)&=&
-g^{m\bar \gamma} \nabla^{\bar q} \, \hat{R}_{\bar q m}{}^\alpha{}_\beta \, \overline{(\nabla_{\gamma} h \, h^{-1})_{\bar\alpha}{}^{\bar\beta}}
-g^{m\bar \gamma} \, (\nabla_m h\, h^{-1})^\alpha{}_{\beta} \, \overline{ \nabla^{\bar q} \hat{R}_{\bar q \gamma}{}^{
\bar\beta}{}_{\bar \alpha}}
\nonumber\\
(II) &=&
g^{m\bar \gamma}g^{p\bar q} T^{r}{}_{mp} \, R_{\bar q r}{}^\alpha{}_{\beta}\, \overline{(\nabla_{\gamma} h \, h^{-1})_{\bar\alpha}{}^{\bar\beta}}
+ g^{m\bar \gamma} \, (\nabla_m h\, h^{-1})^\alpha{}_{\beta} \, \overline{ g^{p\bar q} T^{r}{}_{m p} R_{\bar q r}{}_{\bar\alpha}{}^{\bar\beta}}
\nonumber\\
(III) &=&
 g^{m\bar \gamma} \, \nabla_m Q^\alpha{}_{\beta} \,\overline{(\nabla_{\gamma} h \, h^{-1})_{\bar\alpha}{}^{\bar\beta}} + g^{m\bar \gamma} \, (\nabla_m h\, h^{-1})^\alpha{}_{\beta} \, \overline{\nabla_m Q_{\bar\alpha}{}^{\bar\beta} }
\nonumber\\
(IV)&=&
 -g^{m\bar \gamma} \, \nabla_m \Phi^\alpha{}_{\beta} \,\overline{(\nabla_{\gamma} h \, h^{-1})_{\bar\alpha}{}^{\bar\beta}} - g^{m\bar \gamma} \, (\nabla_m h\, h^{-1})^\alpha{}_{\beta} \, \overline{\nabla_m  \Phi_{\bar\alpha}{}^{\bar\beta} }
\nonumber\\
(V) &=&
 -(\nabla_m h \, h^{-1} )^\beta{}_{\ell}\, \overline{(\nabla_\gamma h \, h^{-1} )^\mu{}_{\alpha}} \left\{Q^{m\bar \gamma} g_{\bar \mu\beta} \, g^{\ell \bar \alpha} - g^{m\bar \gamma} \, Q_{\bar \mu \beta}\, g^{\ell\bar \alpha} + g^{m\bar \gamma} \, g_{\bar \mu \beta} \, Q^{\ell\bar \alpha} \right\}
 \nonumber\\
(VI) &=&
- (\nabla_m h \, h^{-1} )^\beta{}_{\ell}\, \overline{(\nabla_\gamma h \, h^{-1} )^\mu{}_{\alpha}} \left\{- \Phi^{m\bar \gamma} g_{\bar \mu\beta} \, g^{\ell \bar \alpha} + g^{m\bar \gamma} \,  \Phi_{\bar \mu \beta}\, g^{\ell\bar \alpha} - g^{m\bar \gamma} \, g_{\bar \mu \beta} \, \Phi^{\ell\bar \alpha} \right\}.
 \nonumber
 \eea

\subsubsection{Estimate of $S$}
In this section, we estimate terms in the previous expression. The final result will be
\be \label{S-est1}
(\p_t - \Delta) S \leq - {1\over 2} \left(|\bar\nabla(\nabla h \, h^{-1} )|^2+|\nabla(\nabla h \, h^{-1} )|^2 \right) + C(1+S+|T|^2S),
\ee
where $C$ depends on $(X,\hat{g})$, $\Phi$, and bounds for the metric $g$ above and below in terms of the reference $\hat{g}$. To prove this, we will need to estimate terms one by one.
\begin{itemize}
\item[] \underline{For (I):} Because of the presence of connection in e.g. $\nabla^{\bar q} \, \hat{R}_{\bar q m}{}^\alpha{}_\beta$, the first covariant derivatives are of the order $O(S^{1\over 2})$. Therefore,
\bea
(I) \leq C_1 \, (S+S^{1/2}).
\eea

\item[] \underline{For (II):} Recall that
\bea
\hat{R}_{\bar q m}{}^\alpha{}_\beta = R_{\bar q m}{}^\alpha{}_\beta +  \bar \nabla_{q} \left(\nabla_m h \, h^{-1} \right)^\alpha{}_\beta.
\eea
Using this, we can estimate (II) as
\bea
(II) \leq C_2 \, S^{1/2} |T|  + C_2 \, S^{1/2} |T| \, \left\{ |\bar\nabla (\nabla h\, h^{-1})|+ |\bar\nabla (\nabla h\, h^{-1})| \right\}.
\eea

\item[] \underline{For (III):} We have
\be
\nabla_m Q^\alpha{}_\beta = -{1 \over 2} \nabla_m T^\alpha{}_{pq} \bar{T}_\beta{}^{pq} -{1 \over 2} T^\alpha{}_{pq} \nabla_m \bar{T}_\beta{}^{pq}.
\ee
By definition
\be
\nabla_m T^\alpha{}_{pq} = \nabla_m (\Gamma^\alpha_{pq} - \Gamma^\alpha_{qp}).
\ee
By the relation between the reference and evolving connections (\ref{diff-connections2}), we have
\be
\nabla_m T^\alpha{}_{pq} = \nabla_m \hat{T}^\alpha_{pq} + \nabla_m  ( ( \nabla_p h h^{-1})^\alpha{}_q -  ( \nabla_q h h^{-1})^\alpha{}_p).
\ee
Terms $\nabla_m \bar{T}_\beta{}^{pq}= \overline{g^{i \bar{p}} g^{j \bar{q}} \nabla_{\bar{m}} T_{\bar{\beta} ij}}$ can be analyzed similarly. Altogether, we have
\bea
(III) \leq C_3\, (S+S^{1/2}) |T| + C_3 \, S^{1/2} |T| \,  \left\{ |\nabla (\nabla h\, h^{-1})|+ |\bar\nabla (\nabla h\, h^{-1})| \right\}.
\eea

\item[] \underline{For (IV):} Similarly as (I), we have
\bea
(IV) \leq C_4 \, (S+S^{1/2}),
\eea
where the term $\nabla \Phi(z,g(t))$ contributes an order of $S^{1/2}$.

\item[] \underline{For (V):} We directly estimate
\bea
(V) \leq C_5 \, S |T|^2.
\eea

\item[] \underline{For (VI):} We directly estimate
\bea
(VI) \leq C_6 \, S.
\eea
\end{itemize}

Altogether, we obtain
\bea
(\p_t - \Delta)S &\leq& - |\bar\nabla(\nabla h \, h^{-1} )|^2 -|\nabla(\nabla h \, h^{-1} )|^2 \nonumber\\
&&+C \, S^{1/2} |T| \,  \left\{ |\nabla (\nabla h\, h^{-1})|+ |\bar\nabla (\nabla h\, h^{-1})| \right\} \nonumber\\
&& + C S |T|^2 + C(S+S^{1/2})(1+|T|) .
\eea
From here, we can use $2ab \leq a^2 + b^2$ to obtain (\ref{S-est1}) and hence prove Proposition \ref{prop-S-estimate}.

\section{Improved $C^1$ estimate}
\setcounter{equation}{0}
We now assume that $X$ admits a nowhere vanishing holomorphic $(n,0)$ form $\Omega$. This extra structure will allow us to improve the estimate of the previous section.

\begin{theorem} \label{thm-c1-estimate-v2}
Let $(X,\hat{g})$ be a compact Hermitian manifold with nowhere vanishing holomorphic $(n,0)$ form $\Omega$. Let $g(t)$ evolve by Type IIB flow (\ref{typeiib-flow}). Suppose on $[0,t_0]$, we have the estimate
\be
K_1^{-1} \hat{g} \leq g(t) \leq K_1 \hat{g}, \ \ \ |\nabla \log \| \Omega \|_g|_g \leq K_2
\ee
for $K_1,K_2>0$. Then
\be
| \hat{\nabla} g|^2_g \leq C
\ee
on $[0,t_0]$, where $C$ depends on $K_1$, $K_2$, $(X,\hat{g},\Omega)$, $\Phi$, $g(0)$.
\end{theorem}

This theorem is the $C^1$ estimate stated in Theorem 1. Before giving the proof, we compute the evolution of the dilaton function $\log \| \Omega \|_g$.

\subsection{Evolution of the dilaton}
Recall that locally, $\Omega = \Omega(z) dz^1 \wedge \dots \wedge dz^n$ for a local holomorphic function $\Omega(z)$, and
\be
\| \Omega \|^2_g = {\Omega(z) \overline{\Omega(z)} \over \det g}.
\ee
We compute
\be
\p_t \log \| \Omega \|_g = - {1 \over 2} \p_t \log \det g = -{1 \over 2} g^{p \bar{q}} \dot{g}_{\bar{q} p}.
\ee
By the equation of the flow (\ref{typeiib-flow}), we have
\be \label{evol-dilaton0}
\p_t \log \| \Omega \|_g = {1 \over 2} R + {1 \over 4} |T|^2 + {1 \over 2} i \Lambda_\omega \Phi,
\ee
where $\Lambda_\omega$ is defined in (\ref{Lambda-defn}). On the other hand,
\be
R = - g^{p \bar{q}} \partial_p \partial_{\bar{q}} \log \det g = \Delta \log \| \Omega \|_g^2.
\ee
Therefore
\be \label{evol-dilaton}
(\p_t - \Delta) \log \| \Omega \|_g = {1 \over 4} |T|^2 + {1 \over 2} i \Lambda_\omega \Phi.
\ee

\subsection{Maximum principle}
Consider the test function
\be
G = \log S + \varepsilon {\rm Tr} \, h - A \log \| \Omega \|_g,
\ee
for $\varepsilon,A >0$ to be determined. The evolution of $G$ is given by
\be
(\p_t - \Delta)G = {(\p_t - \Delta) S \over S} + {|\nabla S|^2 \over S^2} + \varepsilon (\p_t - \Delta) {\rm Tr} \, h - A (\p_t - \Delta) \log \| \Omega \|_g.
\ee
By (\ref{Trh-est}), (\ref{S-est1}), (\ref{evol-dilaton}),
\be
(\p_t - \Delta) G \leq \left[ {C \over S} + C + C |T|^2 \right] + {|\nabla S|^2 \over S^2} + \varepsilon \left[ -{1 \over K_1^2} S + C \right] - A \left[ {1 \over 4} |T|^2 + {1 \over 2} i \Lambda_\omega \Phi \right].
\ee
Suppose $G$ attains a maximum at a point $(\hat{x},\hat{t})$ with $\hat{t}>0$ and $S(\hat{x},\hat{t})>1$. Then
\be \label{main-ineq1}
0 \leq {|\nabla S|^2 \over S^2} - {\varepsilon \over K_1^2} S -\left({A \over 4} - C \right)|T|^2 + C(\varepsilon,A)
\ee
at $(\hat{x},\hat{t})$. We also have the critical equation $\nabla G (\hat{x},\hat{t})=0$ which implies the relation
\be
{\nabla S \over S} = -\varepsilon \nabla {\rm Tr} \, h + A \nabla \log \| \Omega \|_g.
\ee
Therefore
\be
{|\nabla S|^2 \over S^2} \leq 2 \varepsilon^2  |\nabla {\rm Tr} \, h|^2 + 2 A^2 |\nabla \log \| \Omega \|_g|^2.
\ee
Since
\be
|\nabla {\rm Tr} \, h|^2 = | {\rm Tr} \, \hat{\nabla} h|^2 = g^{i \bar{j}} \hat{g}^{p \bar{q}} \hat{g}^{s \bar{r}} \hat{\nabla}_i g_{\bar{q} p} \hat{\nabla}_{\bar{j}} g_{\bar{r} s} \leq n K_1^2 |\hat{\nabla} g|^2,
\ee
we conclude
\be
{|\nabla S|^2 \over S^2} \leq 2n K_1^2 \varepsilon^2 S + 2 A^2 K_2^2.
\ee
Substituting this inequality in (\ref{main-ineq1}), we obtain
\be
0 \leq - \left({\varepsilon \over K_1^2} - 2n K_1^2\varepsilon^2 \right) S -\left({A \over 4} - C_0 \right)|T|^2 + C(\varepsilon,A)+ 2 A^2 K_2^2.
\ee
Let $\varepsilon = 1/(4nK_1^4)$ and $A = 4 C_0$. Then
\be
0 \leq - {1 \over 8n K_1^6}S + C(\varepsilon,A) + 2 A^2 K_2^2.
\ee
It follows that
\be
S \leq C(K_1,K_2),
\ee
at $(\hat{x},\hat{t})$, and hence $S$ is bounded uniformly at all points $(x,t) \in X \times [0,t_0]$. This completes the proof of Theorem \ref{thm-c1-estimate-v2}.

\section{Higher Order Estimates}
\setcounter{equation}{0}

We now establish higher order estimates for the Type IIB flow, using the theory of quasi-linear parabolic systems.

\subsection{Higher order estimates for quasi-linear parabolic system}

We recall the following theorem about the higher order estimates for general linear and quasi-linear parabolic system from the book by Ladyzenskaja et al. (\cite{LU} Theorem 5.1 in Chapter 7, page 586). Consider the following quasi-linear parabolic systems of the form
\bea\label{system}
\u_t = a^{ij}(x, t, \u) \, \u_{i\, j} + \a(x, t, \u, \u_x),
\eea
in which $\u(x, t)= \left( u^1(x, t), \, u^2(x, t), \, \cdots, \, u^N(x, t) \right)$ is an unknown vector function defined in $Q_T = \Omega \times (0, T)$ with $\Omega \subset \R^n$ and $\a(x, t, \u, {\bf p})$ is a given $N$-dimensional vector-valued function with components $a^\ell\left(x, t, \u(x, t), \u_x(x, t)\right)$; $a^{ij}(x, t, \u)$ is an $n\times n$ matrix function satisfying
\bea
\lambda \cdot |\xi|^2 \leq a^{ij} \xi_i \xi_j \leq \Lambda \cdot |\xi|^2
\eea
for any real $\xi = \left( \xi_1, \cdots \xi_n\right)$.

The theorem states the following

\begin{theorem} Let $\u(x, t)\in C^{2, 1}\left( \bar Q_T\right)$ be a solution of the quasi-linear system (\ref{system}). Suppose that
\begin{itemize}
 \item[1)] The functions $a^{ij}(x, t, \u)$ and their derivatives with respect to the $x_k$ and $u^\ell$ are all continuous in the domain
 \bea
 {\mathcal D}= \{(x, t)\in \bar Q_T : \, |\u|\leq M_0, \, |{\bf p}|\leq M_1\}
 \eea
 where $M_0= \max_{Q_T} | \u(x, t)|$ and $M_1= \max_{Q_T} | \u_x(x, t)|$;

 \item[2)] The functions $a^\ell(x, t, \u, {\bf p})$ are continuous in ${\mathcal D}$.
\end{itemize}
Then for any $Q' \subset Q_T$, there exist two positive constants $C$ and $\alpha$ such that
\bea
 \|\u\|_{C^{1+\alpha, {\alpha/2}}(Q')} < C.
\eea
Here $\alpha$ and $C$ depend on $M_0, M_1, \lambda$, the distance from $Q'$ to $\partial \Omega \times [0,T]$, $\|\u(x, 0)\|_{C^{1+\alpha}(\Omega)}$, and the moduli of the continuity in 1) and 2).

\medskip

If, in addition, the functions  $a^\ell(x, t, \u, {\bf p})$ and $a^{ij}(x, t, \u)$ together with its derivatives satisfy a H\"older condition in ${\mathcal D}$ in the arguments $x, t, \u, {\bf p}$ with exponents $\beta, \beta/2, \beta, \beta$ respectively, then
\bea
 \|\u\|_{C^{2+\beta, {1+\beta/2}}(Q')} < C
\eea
for some constant $C$ depending on $M_0, M_1, \lambda$, ${\rm dist}(Q', \partial \Omega \times [0,T])$, the H\"older continuity of $a^{ij}$ and $a^\ell$, and $\|\u(x, 0)\|_{C^{2+\beta}(\Omega)}$.
\end{theorem}

\

The requirement about the derivatives of $a^{ij}(x, t, \u)$ in condition 1) is used to re-write the system of equations (\ref{system}) into a quasi-linear system of divergence form.

\

\subsection{Higher order estimates for the Type II B flow}

We will use the above estimates for general quasi-linear parabolic systems to derive the higher order estimates for our flow (\ref{typeiib-flow})
\be
\p_t g_{\bar k j} = - \tilde{R}_{\bar k j} - {1\over 2} T_{\bar k p q} \, \bar{T}_{j}{}^{pq} - \Phi(z,g(t))_{\bar k j}.
\ee
From the definition, we know
\bea
\tilde{R}_{\bar k j} =- g^{q\bar p} \, \p_{\bar p} \p_q \, g_{\bar k j} + g^{q\bar p} \, g^{n\bar \ell}\, \p_{\bar p} g_{\bar k n} \, \p_q g_{\bar \ell j}
\eea
and
\bea
T_{\bar k j m} = \p_j g_{\bar k m} - \p_m g_{\bar k j}.
\eea

Therefore, we can treat flow (\ref{typeiib-flow}) as a quasi-linear parabolic system of the form (\ref{system}) by taking $\u = g$,
\bea
a^{ij} (x, t, \u) = g^{-1}
\eea
 and
 \bea
 \a(x, t, \u, {\bf p}) =  g^{q\bar p} \, g^{n\bar \ell}\, \p_{\bar p} g_{\bar k n} \, \p_q g_{\bar \ell j}- {1\over 2} T_{\bar k p q} \, \bar{T}_{j}{}^{pq} - \Phi(z,g(t))_{\bar k j}
 \eea
Using the $C^1$ estimate on $g$, we can check that $a^{ij}$ and $\a$ satisfy the conditions 1) and 2) in the general theorem. Then, the higher order estimates for $g$ follows directly.

\begin{theorem} \label{thm-higher-estimate}
Let $(X,\hat{g})$ be a compact Hermitian manifold with nowhere vanishing holomorphic $(n,0)$ form $\Omega$. Let $g(t)$ evolve by Type IIB flow (\ref{typeiib-flow}). Suppose on $[0,t_0]$, we have the estimate
\be
K^{-1} \hat{g} \leq g(t) \leq K \hat{g}, \ \ \ |\nabla \log \| \Omega \|_g| \leq K
\ee
for $K>0$. Then
\be
\| g \|_{C^{2+\alpha, 1+ \alpha/2}(X \times [0,t_0])} \leq C,
\ee
where $C$ depends on $K$, $(X,\hat{g},\Omega)$, $\Phi$, $g(0)$.
\end{theorem}

Higher order estimates on $g$ follow from the parabolic Schauder estimates. Once $g(t)$ is uniformly bounded in all $C^k$ norms on $X \times [0,t_0]$, a standard compactness argument using the Arzela-Ascoli theorem and the short-time existence theorem gives the extension of the flow to $[0,t_0+\epsilon]$ for some $\epsilon>0$.

\section{Proof of Theorem \ref{thm-zero-source}}
\setcounter{equation}{0}

In comparison with the estimates in Theorem \ref{thm-long-time}, the assumption on the lower bound of the evolving metric $\omega(t)$ is removed when $\Psi=0$.  Indeed, it was noted in \cite{FePi} that the following bound holds along the Anomaly flow:
\be \label{dilaton-upper}
\| \Omega \|_{\omega(t)} \leq \sup_X \| \Omega \|_{\omega_0}.
\ee
This can also be seen from (\ref{evol-dilaton}), which shows that
\be
(\p_t - \Delta) \log \| \Omega \|_{\eta(t)} \geq 0,
\ee
where $\eta = \| \Omega \|_\omega \omega$. The bound (\ref{dilaton-upper}) implies
\be
{\det \, g_0 \over \det \, g} \| \Omega \|^2_{g_0} \leq C.
\ee
If we let $h^i{}_j = (g_0)^{i \bar{k}} g_{\bar{k} j}(t)$, this implies
\be
\det \, h \geq C_0^{-1}.
\ee
If $0<\lambda_n \leq \dots < \lambda_1$ are the eigenvalues of $h$, by assumption we have that
\be
\lambda_1 \leq K_1.
\ee
It follows that
\be
{1 \over \lambda_n} = {\lambda_{n-1} \cdots \lambda_1 \over \det \, h} \leq C_0 K_1^{n-1}.
\ee
Therefore
\be
\omega(t) \geq C^{-1} \omega_0,
\ee
and we may apply Theorem \ref{thm-long-time}. Q.E.D.

\newpage

Department of Mathematics $\&$ Computer Science, Rutgers, Newark, NJ 07102, USA

\smallskip

teng.fei@rutgers.edu

\bigskip

Department of Mathematics, Columbia University, New York, NY 10027, USA

\smallskip

phong@math.columbia.edu

\bigskip

Department of Mathematics, Harvard University, Cambridge, MA 02138, USA

\smallskip

spicard@math.harvard.edu

\bigskip

Department of Mathematics, University of California, Irvine, CA 92697, USA

\smallskip
xiangwen@math.uci.edu

\end{document}